\documentclass[11]{article}
\usepackage[margin=1in]{geometry}
\usepackage{amsfonts,amsmath,amssymb}
\numberwithin{equation}{section}

\title{Trotter-Kato Theorem for Bi-Continuous Semigroups and Approximation of PDEs}
\author{Abdulhameed Qahtan Abbood Altai \\ \small
Department of Mathematics and Computer Science, University of Babylon \\
\small Babil 51002, Iraq \\ \small ahbabil1983@gmail.com}
 
\begin{document}
\maketitle

\paragraph{Abstract}
\
 
In this paper, we introduce formulations of the Trotter-Kato theorem for approximation of bi-continuous semigroups that provide a useful framework whenever convergence of numerical approximations to solutions of PDEs are studied  with respect to an additional locally convex topology coarser than the norm topology to treat the lack of the strong continuity. Applicability of our results is demonstrated using a heat equation with infinite boundaries.


\section{Introduction} \ \ \ \
The theory of strongly continuous semigroups on Banach spaces was created by E. Hille [4], [5] and K. Yosida [10] in order to treat initial value problems for partial differential equations. 

The semigroups of bounded linear operators on a Banach space which are locally bi-equicontinuous with respect to an additional locally convex topology $\tau$ coarser than the $\Vert \cdot \Vert$-topology and the orbit maps $t \rightarrow T (t)x$ are also continuous with respect to this topology  are called bi-continuous semigroups. This concept was proposed by Kühnemund to treat semigroups that are not strongly continuous on a Banach space. She also showed that these kind of semigroups permit systematic theory like Hille-Yosida and Trotter-Kato theorems [7]. 
 
 Albanese and Mangino improved the results of Kühnemund to completely characterize the convergence of bi-continuous semigroups with respect to the topology $\tau$. They obtained a Lie-Trotter product formula and applied it to second order differential operators on $C_{b}(\mathbb{R}^{n})$ [1].

 Ito and Kappel proposed a version of the Trotter-Kato theorem for approximating a linear $C_{o}$-semigroup on a Banach space, which are useful for studying convergence of numerical approximations of solutions to partial differential equations [6]. 
 
 In this paper, we generalize the results of Ito and Kappel and follow the approach of Kühnemund and Albanese and Mangino to propose a version of the Trotter-Kato theorem for approximating a bi-semigroup $T(t)$ on a Banach space $X$ which are useful for studying $\tau$-convergence of numerical $\tau$-approximations of $\Vert \cdot \Vert$-bounded solutions to partial differential equations.
 
In section 2 we introduce some notation and recall some assumptions, definitions and results from [7],[8] and [6]. In section 3 we present a version of the Trotter-Kato theorem for bi-continuous semigroups. In section 4 we discuss how to apply the stability and $\tau$-consistency of the Trotter-Kato theorem for bi-continuous semigroups. And in section 5 we give an illustrative example of a second order heat equation in one space dimension as an application of our results.

\section{Preliminaries} 
\paragraph{Assumptions 2.1.} Let $(X, \Vert \cdot \Vert)$ be a Banach space with topological dual $X^\prime$ and
let $\tau$ be a locally convex topology on $X$ with the following properties:
\begin{enumerate}
\item The space $(X,\tau)$ is sequentially complete on $\Vert \cdot \Vert$-bounded sets, i.e., every 
$\Vert \cdot \Vert$-bounded $\tau$-Cauchy sequence converges in $(X,\tau)$,
\item The topology $\tau$ is Hausdorff and coarser than the $\Vert \cdot \Vert$-topology,
\item The space $(X,\tau)^\prime$ is norming for $(X, \Vert \cdot \Vert)$, i.e.,
$$ \Vert x \Vert = \mathrm{sup} \left\lbrace \vert \langle x,\phi \rangle \vert: \phi \in (X,\tau)^\prime, \Vert \phi \Vert_{(X,\Vert \cdot \Vert)^\prime} \leq 1 \right\rbrace. $$
\end{enumerate}

Let $L(X)$ be a spsce of bounded linear operators defined on $(X, \Vert \cdot \Vert)$ and $P_{\tau}$ be a family of continuous seminorms $sn$ that induce a locally convex topology $\tau$ on $X$. We assume without restrict of generality that
$sn(x) \leq \Vert x \Vert$ for all $x \in X$ and $ sn \in P_{\tau}$ because the $\tau$-topology is coarser than the $\Vert \cdot \Vert$-topology.
 
\paragraph{Definition 2.2.} An operator family $\{T(t) : t \geq 0 \} \subseteq L(X)$ is called (globally)
bi-equicontinuous if for every $\Vert \cdot \Vert$-bounded sequence $(x_{n})_{n \in \mathbb{N}} \subseteq X$ which is $\tau$-convergent to $x \in X$ then 
$$\tau - \lim \limits_{n \to \infty} (T(t)(x_{n}-x)) = 0 $$
uniformly for all $t \geq 0$. It is called locally bi-equicontinuous if for every $t_{0} \geq 0$ the subset $\{ T(t) : 0 \leq t \leq t_{0} \}$ 
is bi-equicontinuous.
 
\paragraph{Definition 2.3.} An operator family $\{T(t) : t \geq 0 \} \subseteq L(X)$ is called a bi-continuous semigroup (with respect to $\tau$ and of type $\omega$) if the following conditions hold
\begin{enumerate}
\item $T(0) = I$ and $T(t + s) = T(t)T(s)$ for all $s, t \geq 0$,
\item The operators $T(t)$ are exponentially bounded, i.e., $\Vert T(t) \Vert_{L(X)} \leq M e^{\omega t}$
for all $t \geq 0$ and some constants $M \geq 1$ and $\omega \in \mathbb{R}$,
\item $(T(t))_{t \geq 0}$ is strongly $\tau$-continuous, i.e., for every $x \in X, \tau-\lim \limits_{t \to 0^{+}} T(t)x = x$, 
\item $(T(t))_{t \geq 0}$ is locally bi-equicontinuous.
\end{enumerate}
 
\paragraph{Definition 2.4.} The $\tau$-generator $A : D(A) \subseteq X \longrightarrow X$ of a bi-continuous semigroup $(T(t))_{t \geq 0}$ on $X$ is the unique operator on $X$ such that its resolvent $R(\lambda,A)$ is
$$R(\lambda,A)x = \int_{0}^{\infty} e^{-\lambda t} T(t) x dt$$
for all $\lambda \in \Lambda_{\omega_{0}} = \{ \lambda \in \mathbb{C} : Re \lambda > \omega_{0} \}$ and $x \in X$. Moreover, the $\tau$-generator $(A ,D(A))$ of a bi-continuous semigroup $(T(t))_{t \geq 0}$ can be defined by
\begin{align*}
Ax = \tau-\lim \limits_{t \to 0^{+}} \frac{T(t)x-x}{t}, \forall x \in D(A)
\end{align*}
where $$ D(A)= \left\lbrace  x \in X: \sup_{(0,1]} \left\Vert  \frac{T(t)x-x}{t} \right\Vert < +\infty \ \mathrm{and} \ \tau-\lim \limits_{t \to 0^{+}} \frac{T(t)x-x}{t} \ \mathrm{exists \ in} \ X \right\rbrace. $$

\paragraph{Definition 2.5.}
Let $\{ (T_{n}(t))_{t \geq 0} : n \in \mathbb{N} \}$ be a sequence of bi-continuous semigroups on $X$. They are
called uniformly bi-continuous if the following conditions hold:
\begin{enumerate}
\item[(i)] There exist $\omega \in \mathbb{R}$ and $M \geq 1$ such that $\Vert T_{n}(t) \Vert \leq M e^{\omega t}$ for all $t \geq 0$ and for all $n \in \mathbb{N}$,
\item[(ii)] $(T_{n}(t))_{t \geq 0}$ are locally bi-equicontinuous uniformly for $n \in \mathbb{N}$, i.e., for every $\Vert \cdot \Vert$-bounded sequence $(x_{k})_{k \in \mathbb{N}} \subset X$ $\tau$-convergent to $x \in X$, that 
\begin{align*}
\tau-\lim \limits_{k \to \infty} {T_{n}(t)(x_{k} - x)} = 0
\end{align*}
uniformly for $t \in [0,t_{o}]$ for every $t_{o} > $ and $n \in \mathbb{N}$.
\end{enumerate}

\paragraph{Definition 2.6.}
A subset $D \subseteq D(A)$ is called a bi-dense if for all $x \in X$ there exists a $\Vert \cdot \Vert$-bounded sequence $(x_{n})_{n \in \mathbb{N}} \subseteq D$ which is $\tau$-convergent to $x$. And it is called a bi-core for $A$ if for all $x \in D(A)$ there exists a sequence $(x_{n})_{n \in \mathbb{N}} \subseteq D$ such that $(x_{n})_{n \in \mathbb{N}}$ and $(Ax_{n})_{n \in \mathbb{N}}$ are $\Vert \cdot \Vert$-bounded and $\tau-\lim \limits_{n \to \infty} x_{n} = x$ and $\tau-\lim \limits_{n \to \infty} Ax_{n} = Ax$. 

\paragraph{Remark 2.8.} If $\rho(A) \neq \emptyset$, then $D$ is a bi-core for $A$ if and only if $D$ and $(\lambda - A)(D)$ are bi-dense in $X$.

\paragraph{Proposition 2.9.} Let $(A,D(A))$ be the $\tau$-generator of a bi-continuous semigroup $(T(t))_{t \geq 0}$
of type $\omega$ on $X$. Then we have
\begin{align*}
\dfrac{d^{k}}{d \lambda^{k}} R(\lambda,A)x = (-1)^{k} \int_{0}^{\infty} t^{k} e^{-\lambda t} T(t) x dt
\end{align*}
for all $x \in X, k \in \mathbb{N}$ and $\lambda \in \Lambda_{\omega_{0}}$. And there exists for each $\omega > \omega_{0}$ a
constant $M \geq 1$ such that
\begin{align*}
\Vert R(\lambda,A)^{k} \Vert \leq \frac{M}{(Re \lambda - \omega)^{k}}
\end{align*}
for all $k \in \mathbb{N}$ and $\lambda \in \Lambda_{\omega}$.

\paragraph{Proposition 2.10.} Let $(T(t))_{t \geq 0}$ be a bi-continuous semigroup on $X$ and $(A ,D(A))$ be the $\tau$-generator of $(T(t))_{t \geq 0}$. Then the following properties hold.
\begin{enumerate}
\item If $x \in D(A)$, then $T(t)x \in D(A)$ for all $t \geq 0, T(t)x$ is continuously differentiable
in $t$ with respect to the topology $\tau$, and
\begin{align*}
\frac{d}{dt}T(t)x = A T(t)x = T(t) A x, \forall t \geq 0
\end{align*}
\item An element $x \in X$ belongs to $D(A)$ and $Ax = y$ if and only if
\begin{align*}
T(t)x - x = \int_{0}^{t} T(s) y ds, \ \forall t \geq 0
\end{align*}
\item The operator $(A,D(A))$ is bi-closed, i.e. for all sequences $(x_{n})_{n \in \mathbb{N}} \subseteq D(A)$ with $(x_{n})_{n \in \mathbb{N}}$ and $(Ax_{n})_{n \in \mathbb{N}}$ $\Vert \cdot \Vert-$bounded, $x_{n} \stackrel{\tau}\longrightarrow x \in X$ and $Ax_{n} \stackrel{\tau}\longrightarrow y \in X$ we have $x \in D(A)$ and $Ax = y$.
\end{enumerate}

\section{The Trotter-Kato Theorem} \ \ \ \
Let $(Z,\Vert.\Vert)$ and $(X_{n},\Vert.\Vert_{n}),n=1,2,...,$ be Banach spaces and $X$ be a closed linear subspace of $Z$ with a locally convex topology $\tau$ induced by a family of continuous seminorms $P_{\tau}$. A bi-continuous semigroup $T(\cdot)$ with infinitesimal $\tau$-generator $A$ is defined on $X$. Our mission is to construct approximationg $\tau$-generators $A_{n}$ on the spaces $X_{n}$ such that the bi-continuous semigroups $T_{n}(\cdot)$ generated by $A_{n}$ $\tau$-approximate $T(\cdot)$ as it will be explained below.

Before we state the Trotter-Kato theorem for approximation of bi-continuous semigroups we need to consider the following:

\paragraph{Assumptions 3.1.} 
For every $n = 1,2,...,$ there exist bounded linear operators $P_{n}:Z \to X_{n}$ and $E_{n}:X_{n} \to Z$ satisfying
\begin{enumerate}
\item[(A1)] $\Vert P_{n} \Vert \leq M_{1}, \Vert E_{n} \Vert \leq M_{2}$, where $M_{1},M_{2}$ are independent of $n$,
\item[(A2)] $(E_{n}P_{n}x)_{n \in \mathbb{N} }$ is $\Vert \cdot \Vert$-bounded and $\tau - \lim \limits_{n \to \infty} (E_{n}P_{n}x-x) = 0 $ for all $x \in X$,
\item[(A3)] $P_{n}E_{n} = I_{n}$, where $I_{n}$ is the identity operator on $X_{n}$.
\end{enumerate}

\paragraph{Lemma 3.2.} 
The mappings $\pi_{n}:Z \to Z_{n}$ defined by $\pi_{n} = E_{n}P_{n},n=1,2,...,$ are projections, i.e., $\pi^{2}_{n} = \pi_{n}$ and $\mathrm{range} \ \pi_{n} = Z_{n} $ where $Z_{n} = \mathrm{range} \ E_{n}$ are closed subspaces of $Z$.

\paragraph{Lemma 3.3.} 
$\tilde{T}_{n}(t) = E_{n}T_{n}(t)P_{n} \vert_{Z_{n}}$ defines a bi-continuous semigroup on $Z_{n}$ with infinitesimal $\tau$-generator $\tilde{A}_{n}$ given by $D (\tilde{A}_{n}) = E_{n} D(A_{n})$ and $\tilde{A}_{n}(t) = E_{n}A_{n}(t)P_{n} \vert_{Z_{n}}$.\

We denote by $A \in G(M,\omega,X)$ for some $M \geq 1$ and $\omega \in \mathbb{R}$ if $A$ is the infinitesimal $\tau$-generator of a bi-continuous semigroup $T(t),t \geq 0$ satisfying $\Vert T(t) \Vert \leq M e^{\omega t}, t \geq 0$. 

\paragraph{Theorem 3.4 (Trotter-Kato).} Assume that $(A1)$ and $(A3)$ are satisfied. Let $T(t)$ and $T_{n}(t)$ be bi-continuous semigroups with $A \in G(M,\omega,X)$ and $A_{n} \in G(M,\omega,X_{n})$ on $X$ and $X_{n}$ respectively, and let $D$ be a bi-core for $A$. Then the following are equivalent:
\begin{enumerate}
\item[(a)] There exists a $\lambda_{o} \in \rho (A) \cap \bigcap^{\infty}_{n = 1} \rho (A_{n})$ such that, for all $x \in X$, $(E_{n} (\lambda_{o} I_{n} - A_{n})^{-1} P_{n}x)_{n \in \mathbb{N}}$ is $\Vert \cdot \Vert$-bounded and
\begin{align*}
\tau- \lim \limits_{n \to \infty} E_{n} (\lambda_{o} I_{n} - A_{n})^{-1} P_{n}x = (\lambda_{o} I - A)^{-1}x
\end{align*}
\item[(b)] For every $x \in X$ and $t \geq 0$, $(E_{n} T_{n} P_{n} x)_{n \in \mathbb{N}}$ is $\Vert \cdot \Vert$-bounded and
$$ \tau- \lim \limits_{n \to \infty} E_{n} T_{n} P_{n} x = T(t)x $$
uniformly on bounded $t$-intervals.
\end{enumerate}

\paragraph{Proof.}
Let $sn \in P_{\tau}$ and set $Z_{n} = \mathrm{range} E_{n}$ and $\pi_{n} = E_{n}P_{n},n = 1,2,...$ . We need to establish the following two statements which are equivalent to $(a)$ and $(b)$ to prove our theorem:
\begin{enumerate}
\item[$(\tilde{a})$] There exists a $\lambda_{o} \in \rho (A) \cap \bigcap^{\infty}_{n = 1} \rho (\tilde{A}_{n})$ such that, for all $x \in X$, $((\lambda_{o} \tilde{I}_{n} - \tilde{A}_{n})^{-1} \pi_{n}x)_{n \in \mathbb{N}}$ is $\Vert \cdot \Vert$-bounded and 
$$ \tau- \lim \limits_{n \to \infty} (\lambda_{o} \tilde{I}_{n} - \tilde{A}_{n})^{-1} \pi_{n}x = (\lambda_{o} I - A)^{-1}x$$
\item[$(\tilde{b})$] For every $x \in X$ and $t \geq 0$, that $(\tilde{T}_{n} \pi_{n} x)_{n \in \mathbb{N}}$ is $\Vert \cdot \Vert$-bounded and
$$ \tau- \lim \limits_{n \to \infty} \tilde{T}_{n} \pi_{n} x = T(t)x $$
uniformly on bounded $t-$intervals.
\end{enumerate}
For convenience, we will write $T_{n}(t)$ and $A_{n}$ instead of $\tilde{T}_{n}(t)$ and $\tilde{A}_{n}$ respectively if there is no confusion for the rest of the proof, and without loss of generality we assume that $(\tilde{a})$ holds for $\lambda_{o} = 0$. \\
i) To show that $(\tilde{a})$ implies $(\tilde{b})$ at first. For $x \in X$, let
$$ e_{n}(t) = (T_{n}(t) \pi_{n} - \pi_{n} T(t))x, n = 1,2,..., t \geq 0.$$
For $x \in D(A)$, the function 
$$u_{n}(t) = A^{-1}_{n} e_{n}(t), n = 1,2,..., t \geq 0 $$
is in $C^{1}(0,\infty;Z_{n})$ and satisfied 
\begin{equation}
\begin{aligned}
\dot{u}_{n} &= A_{n} u_{n} + \pi_{n} \Delta_{n} A T(t)x, \\
u_{n}(0) &= 0,
\end{aligned}		   
\end{equation}
where $\Delta_{n} = A^{-1} - A^{-1}_{n} \pi_{n}$. By proposition 2.10 that $A^{-1}_{n} T_{n}(t) \pi_{n} x =  T_{n}(t) A^{-1}_{n} \pi_{n} x$ is continuously differentiable on $[0,\infty)$ because $A^{-1}_{n} \pi_{n} x$ is in $D(A)$, and $A^{-1}_{n} \pi_{n} T(t) x$ is continuously differentiable because $x \in D(A)$ and $A^{-1}_{n} \pi_{n} : X \to Z_{n}$ is bounded operator.

Using the variation of parameter formula we obtain from $(3.1)$, for $t \geq 0, x \in D(A)$, that
\begin{align}
u_{n}(t) = \int^{t}_{0} T_{n}(t - \tau) \pi_{n} \Delta_{n} A T(\tau)x d\tau.
\end{align}
For $x \in D(A^{2})$, we integrate $(3.2)$ by parts to get
\begin{equation}
\begin{aligned} 
u_{n}(t) = &- A^{-1}_{n} \pi_{n} \Delta_{n} A T(t)x +  A^{-1}_{n} T_{n}(t) \pi_{n} \Delta_{n} A x \\
		   &+ A^{-1}_{n} \int^{t}_{0} T_{n}(t - \tau) \pi_{n} \Delta_{n} A^{2}       
		      T(\tau)x d\tau, t \geq 0.
\end{aligned}		   
\end{equation}
So, the error representation becomes
\begin{equation}
\begin{aligned}
e_{n}(t) = &- \pi_{n} \Delta_{n} A T(t)x + T_{n}(t) \pi_{n} \Delta_{n} A x \\
		   &+ \int^{t}_{0} T_{n}(t - \tau) \pi_{n} \Delta_{n} A^{2} T(\tau)x d\tau, t \geq 0, x \in D(A^{2}).
\end{aligned}		   
\end{equation}

To show that $(e_{n}(t))_{n \in \mathbb{N}}$ is $\Vert \cdot \Vert$-bounded and $\tau- \lim \limits_{n \to \infty} e_{n}(t) = 0$ uniformly on $0 \leq t \leq t_{o}$, we need to consider the terms on the right-hand side of $(3.4)$. Fix  $t_{o} > 0$, the set $\{ T(t)Ax: 0 \leq t \leq t_{o} \}$ is compact and by $(\tilde{a})$ that $(- \pi_{n} \Delta_{n} A T(t)x )_{n \in \mathbb{N}}$ is $\Vert \cdot \Vert$-bounded and 
\begin{align*}
\tau- \lim \limits_{n \to \infty} (- \pi_{n} \Delta_{n} A T(t)x ) = 0
\end{align*}
uniformly on $[0,t_{o}]$. Also by $(\tilde{a})$, that $(T_{n}(t) \pi_{n} \Delta_{n} A x)_{n \in \mathbb{N}}$ is $\Vert \cdot \Vert$-bounded and 
\begin{align*}
\tau- \lim \limits_{n \to \infty}  T_{n}(t) \pi_{n} \Delta_{n} A x = 0
\end{align*}
uniformly on $[0,t_{o}]$ because $ \Vert T_{n}(t) \Vert \leq M e^{\omega t}, t \geq 0, n = 1,2,...$ .

Furthermore, since, for $x \in D(A^{2})$, the set $\{ A^{2} T(t)x: 0 \leq t \leq t_{o} \}$ is compact and by $(\tilde{a})$ that $(\Delta_{n} A^{2} T(t) x)_{n \in \mathbb{N}}$ is $\Vert \cdot \Vert$-bounded and
\begin{align*}
\tau- \lim \limits_{n \to \infty} \Delta_{n} A^{2} T(t) x = 0 
\end{align*}
uniformly on $[0,t_{o}]$, then $(T_{n}(t - \tau) \pi_{n} \Delta_{n} A^{2} T(\tau)x) _{n \in \mathbb{N}}$ is $\Vert \cdot \Vert$-bounded and
\begin{align*}
\tau- \lim \limits_{n \to \infty} \int^{t}_{0} T_{n}(t - \tau) \pi_{n} \Delta_{n} A^{2} T(\tau)x d\tau = 0
\end{align*}
uniformly on $[0,t_{o}]$ by the Lebesgue's dominated convergence theorem. This shows that $(e_{n}(t))_{n \in \mathbb{N}}$ is $\Vert \cdot \Vert$-bounded and $\tau- \lim \limits_{n \to \infty} e_{n}(t) = 0$ uniformly on $[0,t_{o}]$ for any $x \in D(A^2)$. 

Since $\{ T(t)x: 0 \leq t \leq t_{o} \}$ is compact,
\begin{align}
\pi_{n}T(t)x - T(t)x = (\pi_{n} - I) \Delta_{n} AT(t)x,
\end{align}
and by $(\tilde{a})$ that $(\pi_{n} T(t) x)_{n \in \mathbb{N}}$ is $\Vert \cdot \Vert$-bounded and
\begin{align*}
\tau- \lim \limits_{n \to \infty} \pi_{n}T(t)x = T(t)x. 
\end{align*}
uniformly on $[0,t_{o}]$, then $(T_{n}(t) \pi_{n} x)_{n \in \mathbb{N}}$ is $\Vert \cdot \Vert$-bounded and
\begin{align*} 
sn(T_{n}(t) \pi_{n} x - T(t) x) \leq sn(e_{n}(t)) + sn(\pi_{n} T(t) x - T(t) x) \to 0 \ \mathrm{as} \ n \to \infty.
\end{align*}
uniformly on $[0,t_{o}]$. \\
ii) And to show that $(\tilde{b})$ implies $(\tilde{a})$. By $(\tilde{b})$, the choice of $\lambda$ and the Lebesgue's dominated convergence theorem, for $Re \lambda > \omega$, that
\begin{align*}
sn( (\lambda I_{n} - A_{n})^{-1} \pi_{n}x - (\lambda I - A)^{-1}x) \leq \int^{\infty}_{0} e^{-Re \lambda t} sn(T_{n}(t) \pi_{n} x - T(t) x) dt \to 0 \ \mathrm{as} \ n \to \infty
\end{align*}
uniformly on $[0,t_{o}]$. $\square$

\paragraph{Remark 3.5.}
According to theorem 3.4, we will call the assumption $A_{n} \in G(M,\omega,X_{n}),n = 1,2,\ldots$ by the stability property of the $\tau$-approximations and the statement $(a)$ by the $\tau$-consistency property of the $\tau$-approximations. In this case, the Trotter-Kato theorem states that, under the assumption of stability of $\tau$-approximations, $\tau$-consistency is equivalent to $\tau$-convergence.

\section{ESTABLISHING OF STABILITY AND $\tau$-CONSISTENCY} \ \ \ \
Applying the Trotter-Kato theorem could face two major difficulties, verification the stability property and computation of the resolvents $(\lambda I_{n} - A_{n})^{-1}$ in verification of the $\tau-$consistency property. And the best way to exceed these difficulties is to use the dissipativity estimates possibly after renorming the spaces $X_{n}$ with uniformly equivalent norms to verify the stability property and avoiding computation of the resolvent operators $(\lambda I_{n} - A_{n})^{-1}$ and direct verification of the $\tau$-consistency property. So we need to consider the following proposition to replace the condition $(a)$ in Theorem 3.5 by an equivalent condition guarantee $\tau$-convergence of the operators $A_{n}$ to $A$.

\paragraph{Proposition 4.1.}
Let the assumptions of Theorem 3.4 be satisfied. Then statement $(a)$ of Theorem 3.4 is equivalent to $(A2)$ and the following two statements:
\begin{enumerate}
\item[(C1)] There exists a subset $D \subseteq D(A)$ such that $D$ and $(\lambda_{o}I - A)D$ are bi-dense sets in $X$ for a $\lambda_{o} > \omega$.
\item[(C2)] For all $u \in D$ there exists a sequence $(\overline{u}_{n})_{n \in \mathbb{N}}, \overline{u}_{n} \in D(A_{n})$, such that $(E_{n} \overline{u}_{n})_{n \in \mathbb{N}}$ and $(E_{n} A_{n} \overline{u}_{n})_{n \in \mathbb{N}}$ are $\Vert \cdot \Vert$-bounded and
\begin{align*}
\tau- \lim \limits_{n \to \infty} E_{n} \overline{u}_{n} = u \ \mathrm{and} \ \tau- \lim \limits_{n \to \infty} E_{n} A_{n} \overline{u}_{n} = Au.
\end{align*}
\end{enumerate}

\paragraph{Proof.} Without loss of generality we assume that $\lambda_{o} = 0$.\\
i) To show that $(a)$ implies $(A_{2})$ and $(C_{1}),(C_{2})$ at first.  We set $D = D(A)$ to get $AD = X$, i.e., $(C_{1})$ is satisfied. $(3.5)$ guarantees that $(a)$ implies $(A_{2})$. For fixed $u \in D(A)$, choose $x \in X$ with $u = A^{-1}x$ and set $\overline{u}_{n} = A^{-1}_{n}P_{n}Au$. By $(a)$ that $(E_{n} \overline{u}_{n})_{n \in \mathbb{N}}$ is $\Vert \cdot \Vert$-bounded and
\begin{align*}
\tau- \lim \limits_{n \to \infty}  (E_{n} \overline{u}_{n} - u) = \tau- \lim \limits_{n \to \infty} (E_{n} A^{-1}_{n}P_{n}x - A^{-1}x) = 0.
\end{align*}
Moreover, by $(A_{2})$ that $(E_{n} A_{n} \overline{u}_{n})_{n \in \mathbb{N}}$ is $\Vert \cdot \Vert$-bounded and  
\begin{align*}
\tau- \lim \limits_{n \to \infty} (E_{n} A_{n} \overline{u}_{n} - Au) = \tau- \lim \limits_{n \to \infty} (E_{n} P_{n}x - x) = 0
\end{align*}
Therefore $(C_{2})$ is also satisfied.\\
ii) And to show that $(A_{2})$ and $(C_{1}),(C_{2})$ imply $(a)$ we use the identity 
\begin{align}
E_{n} A^{-1}_{n} P_{n} - A^{-1} = E_{n} ( A^{-1}_{n} P_{n} A - P_{n})A^{-1} + (E_{n} P_{n} -I)A^{-1}.
\end{align}
For $x \in AD$, choose $u \in D$ with $x = Au$ and set $u_{n} = A^{-1}_{n}P_{n}x = A^{-1}_{n}P_{n} Au$. In addition to, for $u$, take $\overline{u}_{n}$ as in $(C_{2})$. By $(A_{1})$ and $(C_{2})$ we have $(\overline{u}_{n} - P_{n} u)_{n \in \mathbb{N}}$ is $\Vert \cdot \Vert$-bounded and
\begin{align*}
sn_{n} (\overline{u}_{n} - P_{n} u) = sn_{n} (P_{n} (E_{n} \overline{u}_{n} - u)) \leq M_{1} sn ((E_{n} \overline{u}_{n} - u)) \to 0
\end{align*} 
as $n \to \infty$. Also, $(\overline{u}_{n} - u_{n})_{n \in \mathbb{N}}$ is $\Vert \cdot \Vert$-bounded and
\begin{align*}
sn_{n} (\overline{u}_{n} - u_{n}) & \leq sn_{n} (A^{-1}_{n}) sn_{n} (A_{n} \overline{u}_{n} - P_{n} Au) \\
								 & \leq sn_{n} (A^{-1}_{n}) sn_{n} (P_{n}) sn (E_{n} A_{n} \overline{u}_{n} -  Au) \to 0
\end{align*}
as $n \to \infty$. Thus, $(P_{n}u - u_{n})_{n \in \mathbb{N}}$ is $\Vert \cdot \Vert$-bounded and
\begin{align*}
sn_{n} (P_{n}u - u_{n}) \leq sn_{n} (P_{n} u - \overline{u}_{n}) + sn_{n}(\overline{u}_{n} - u_{n}) \to 0
\end{align*}
as $n \to \infty$. The above estimate together with $(4.1)$ and $(A_{2})$ implies that $(E_{n} A^{-1}_{n} P_{n}x)_{n \in \mathbb{N}}$ is $\Vert \cdot \Vert$-bounded and
\begin{align*}
sn (E_{n} A^{-1}_{n} P_{n}x - A^{-1}x) &\leq sn (E_{n}(u_{n} - P_{n}u)) + sn (E_{n}P_{n}u - u) \\
&\leq M_{2} sn_{n} (u_{n} -P_{n}u) + sn (E_{n}P_{n}u - u) \to 0
\end{align*}
as $n \to \infty$ for all $x \in AD$. $\square$

\section{Examples} 

\paragraph{5.1. The heat equation.} Consider the heat equation 
\begin{align}
\begin{array}{ll} 
\dfrac{\partial}{\partial t} u(t,x) + \dfrac{\partial^{2}}{\partial x^{2}} u(t,x) = 0, - \infty < x < \infty,  \\
u(t,x) \rightarrow 0 \ \mathrm{as} \ |x| \rightarrow \infty
\end{array}		   
\end{align}
Let $X = C_{b,0}(\mathbb{R})$ be the space of bounded continuous functions on $\mathbb{R}$ vanishing at $x \rightarrow \infty \ \mathrm{and} \ x \rightarrow -\infty$ endowed with the supremum norm $\Vert\cdot\Vert_{\infty}$. $X$ is a locally convex space with respect to the topology $\tau$ induced by a family $P_{\tau}$ of seminorms defined by
\begin{align*}
sn_{l}(f) = \sup_{ |x| \leq l} |f(x)|,  
\end{align*}
with a locally convex base constructed by finite intersections of the sets,
\begin{align*}
V = \left\lbrace f \in X: sn_{l}(f) < \frac{1}{l} \right\rbrace,  \ l \in \mathbb{N},
\end{align*}
as $\mathbb{R} = \bigcup_{l = 1}^{\infty} [-l,l]$ (see [2],[9]) and satisfies the assumptions 2.1 (see [7],[1]). 
The linear operator $A$ defined by
\begin{align*} 
Au = - u^{\prime \prime}, u \in D(A)
\end{align*}
with $ D(A) = \left\lbrace u \in X | u \ \mathrm{is \ absolutely \ continuous \ on} \ \mathbb{R} \ \mathrm{with} \ u^{\prime \prime} \in X \ \mathrm{and} \ u(t,x) \rightarrow 0 \ \mathrm{as} \ |x| \rightarrow \infty  \right\rbrace $ generates a bi-semigroup on $X$ but not a $C_{0}$-semigroup (see [8],[3]). 
We use the second order finite difference scheme 
\begin{align}
\begin{array}{ll}
\dfrac{d}{dt}u_{k}(t) =  - \dfrac{u_{k+1}(t) - 2u_{k}(t) + u_{k-1}(t)}{{(\bigtriangleup x)}^{2}}, \ k = -l N + 1, \ldots, l N - 1,
\end{array}
\end{align}
where $ \mathrm{col}(u_{-l N + 1},\ldots,u_{l N}) \in X_{n} = \mathbb{R}^{n}$ with $n = card \{ -l N + 1, \ldots, l N \}$ and $u_{k}(t)$ represents an approximating value for $u(t,x)$ at $x_{k} = k \bigtriangleup x,  \bigtriangleup x = 1/N$ which  discretizes $[-l,l]$. From equation $(5.2)$ the approximating generators $A_{n}$ on $\mathbb{R}^{n}$ are given by
\begin{align*}
(A_{n}u)_{k} = - \dfrac{u_{k+1} - 2u_{k} + u_{k-1}}{{(\bigtriangleup x)}^{2}}, \ k = -l N + 1, \ldots, l N - 1.
\end{align*} 
Let $P_{n}, E_{n}$ and $\Vert \cdot \Vert_{n}$  be defined as
\begin{align*}
E_{n} u &= \sum_{k = -lN + 1}^{lN} u_{k} B_{k}(x), \ u \in X_{n}, \\
(P_{n} \phi)_{k} &= \phi (x_{k}), k = -lN + 1,\ldots,lN, \ \phi \in X, \\
\Vert u \Vert_{n} &= \max_{k = -lN + 1,\ldots,lN}  \vert u_{k} \vert,  \ u \in X_{n}, 
\end{align*}
where the first order $B$-spline $B_{k}(x), k = -l N + 1, \ldots, l N,$ for $x \in [-l,l]$ is given by
\begin{align} B_{k}(x) = \left\{ \begin{array}{ll}  
N(x - x_{k}), &\mbox{$x \in [x_{k -1},x_{k}]$}, \\ 
N(x_{k+1} - x),  &\mbox{$x \in [x_{k},x_{k+1}]$}, \\
0              &\mbox{$\mathrm{otherwise}$}.
\end{array}
\right.
\end{align} 
It is easy to satisfy the assumptions $(A1)-(A3)$. The elements $v$ of the duality set $F_{n}(u) \subset X^{*}_{n}$ for $u \in X_{n}$ are given by
\begin{align*} v_{k} = \left\{\begin{array}{ll}
\Vert u \Vert_{n} \mathrm{sgn} \ u_{i}  &\mbox{ $\mathrm{if} k = i$ }, \\ 
0 &\mbox{ $\mathrm{if} k \neq i$ },
\end{array}
\right.
\end{align*}
where $i$ is an index such that $|u_{i}| = \max_{k} |u_{k}|$. Then the stability property is satisfied because $<A_{n}u,v> \leqslant 0$ for all $v \in F_{n}(u)$. For the $\tau$-consistency property, let $D = D(A)$ and $\overline{u}_{n} \in X_{n}$ for $u \in D(A)$  which establishes the condition $(C_{1})$ in Proposition 4.1 with $\omega = 0$. And for the condition $(C_{2})$, since $E_{n} \overline{u}_{n}$ is the first order spline interpolating the continuously differentiable function $u$ at the meshpoints, then, for $ sn \in P_{\tau}$, $(E_{n} \overline{u}_{n})_{n \in \mathbb{N}}$ is $ \Vert \cdot \Vert_{\infty}$-bounded and
\begin{align*}
sn(E_{n} \overline{u}_{n} - u) \rightarrow 0 \ \mathrm{as} \ N \rightarrow \infty.
\end{align*}
Moreover, the sequence  $(E_{n} A_{n} \overline{u}_{n})_{n \in \mathbb{N}}$ is $ \Vert \cdot \Vert_{\infty}$-bounded and, for $\xi_{k} \in (x_{k-1},x_{k}),\zeta_{k} \in (x_{k},x_{k+1})$ and $\eta_{k} \in (\xi_{k},\zeta_{k})$, that
\begin{align*}
sn(E_{n} A_{n} \overline{u}_{n} - A u) &= sn \left( u^{\prime\prime} - \sum_{k = -lN + 1}^{lN - 1} \dfrac{u(x_{k+1}) - 2u(x_{k}) + u(x_{k-1})}{(\bigtriangleup x)^{2}}  B_{k} \right) \\
&= sn \left( u^{\prime\prime} - \sum_{k = -lN + 1}^{lN} \dfrac{u^{\prime} (\zeta_{k}) - u^{\prime} (\xi_{k})}{\bigtriangleup x} B_{k} \right) \\
&= sn \left( u^{\prime\prime} - \sum_{k = -lN + 1}^{lN} u^{\prime\prime} (\eta_{k}) B_{k} \right) \\
&\leq sn \left( u^{\prime\prime} - \sum_{k = -lN + 1}^{lN} u^{\prime\prime} (x_{k}) B_{k} \right) + \max_{k = -lN + 1,\ldots,lN} |u^{\prime\prime} (x_{k}) - u^{\prime\prime} (\eta_{k})| \\
&\leq sn \left( u^{\prime\prime} - \sum_{k = -lN + 1}^{lN} u^{\prime\prime} (x_{k}) B_{k} \right) + \omega (u^{\prime\prime}; (\bigtriangleup x)^{2}) \rightarrow 0 \ \mathrm{as} \ n \rightarrow \infty
\end{align*}
where $\sum^{lN}_{k = -lN + 1} u^{\prime\prime} (x_{k}) B_{k}$ is the first order spline interpolating $u^{\prime\prime}$ at the meshpoints and $ (\bigtriangleup x)^{2} \rightarrow \omega (u^{\prime\prime}; (\bigtriangleup x)^{2})$ denotes the modulus of continuity for $u^{\prime\prime}$.

\end{document}